\documentclass{amsart}         
\usepackage{amscd}             

\renewcommand{\H}{{\mathbb H}}
\renewcommand{\Im}{\operatorname{Im}}
\newcommand{\Isom}{\operatorname{Isom}}

\newcommand{\CS}{\operatorname{CS}}
\newcommand{\tracefield}{\operatorname{TraceField}}

\newcommand{\PGL}{\operatorname{PGL}}
\newcommand{\sgn}{\operatorname{sgn}}
\newcommand{\vol}{\operatorname{vol}}
\newcommand{\Tor}{\operatorname{Tor}}
\renewcommand{\P}{{\mathcal P}}
\newcommand{\B}{{\mathcal B}}
\newcommand{\C}{{\mathbb C}}
\newcommand{\Z}{{\mathbb Z}}
\newcommand{\Q}{{\mathbb Q}}
\newcommand{\R}{{\mathbb R}}
\newcommand{\CP}{{\C\mathbb P}}
\newcommand{\calC}{{\mathcal C}}
\newcommand{\calQ}{{\mathcal Q}}
\newcommand{\calZ}{{\mathcal Z}}
\newcommand{\subdot}{_\bullet}
\newcommand{\bH}[1][n]{\partial\overline{\H}^{#1}}

\newtheorem{theorem}{Theorem}[section]
\newtheorem{lemma}[theorem]{Lemma}
\newtheorem{proposition}[theorem]{Proposition}
\newtheorem{corollary}[theorem]{Corollary}

\theoremstyle{definition}
\newtheorem{definition}[theorem]{Definition}
\newtheorem{example}[theorem]{Example}

\newtheorem{remark}[theorem]{Remark}
\newtheorem*{remark*}{Remark}

\newtheorem{condition}[theorem]{Condition}

\begin{document}
\title[Bloch invariants of hyperbolic 3-manifolds]{Bloch invariants of
hyperbolic 3-manifolds}
\author{Walter D. Neumann}
\address{Department of Mathematics\\The University of
Melbourne\\Parkville, Vic 3052\\Australia}
\email{neumann@maths.mu.oz.au}
\author{Jun Yang}
\address{Department of Mathematics\\Duke University\\Durham NC 27707}
\email{yang@math.duke.edu}
\thanks{This research is supported by the Australian Research Council
and the U.S. National Science Foundation.}

\begin{abstract}
We define an invariant $\beta(M)$ of a finite volume hyperbolic 3-manifold $M$
in the Bloch group $\B(\C)$ and show it is determined by the simplex
parameters of any degree one ideal triangulation of $M$. We show $\beta(M)$
lies in a subgroup of $\B(\C)$ of finite $\Q$-rank determined by the
invariant trace field of $M$.  Moreover, the Chern-Simons invariant of $M$ is
determined modulo rationals by $\beta(M)$.  This leads to a simplicial formula
and rationality results for the Chern Simons invariant which appear
elsewhere.

Generalizations of $\beta(M)$ are also described, as well as several
interesting examples.  An appendix describes a scissors congruence
interpretation of $\B(\C)$.
\end{abstract}
\maketitle
\section{Introduction}

Let $M=\H^3/\Gamma$ be an oriented hyperbolic manifold of finite
volume (so $\Gamma$ is a torsion free Kleinian group).  It is known
that $M$ has a degree one ideal triangulation by ideal simplices
$\Delta_1,\ldots, \Delta_n$ (see sect.\ \ref{ideal triang}).  Let
$z_i\in\C$ be the parameter of the ideal simplex $\Delta_i$ for each
$i$.  These parameters define an element $\beta(M)=\sum_{i=1}^n[z_i]$
in the pre-Bloch group $\P(\C)$ (as defined below and in
\cite{dupont-sah}, for example).
\begin{theorem}\label{theorem1}
The above element\/ $\beta(M)$ can be defined without reference to the
ideal decomposition, so it 
depends only on $M$. Moreover, it lies in the
Bloch group\/ $\B(\C)\subset\P(\C)$.
\end{theorem}

The independence of $\beta(M)$ on ideal triangulation holds even
though our concept of degree one ideal triangulation is rather more
general than ideal triangulation concepts often considered.

We prove this theorem as follows.  There is an exact sequence (mod
$2$-torsion) due to Bloch and Wigner (cf.~\cite{dupont-sah})
$$0\to\mu\to H_3(\PGL(2,\C);\Z)\to \B(\C)\to 0,$$ where
$\mu\subset\C^*$ is the group of roots of unity.  If $M$ is compact
then there is a ``fundamental class'' $[M]\in H_3(\PGL(2,\C);\Z)$ and
we show $\beta(M)$ is the image of $[M]$ in $\B(\C)$.  We do this by
factoring through a certain relative homology group
$H_3(\PGL(2,\C),\CP^1;\Z)$ for which the relationship between $[M]$
and $\beta(M)$ is easier to see (in fact Dupont and Sah
\cite{dupont-sah} show this relative homology group maps
isomorphically to $\P(\C)$).  In the non-compact case we also find a
fundamental class $[M]$ in this relative homology group that maps to
$\beta(M)\in\P(\C)$, thus proving that $\beta(M)$ is independent of
triangulation. The fact that it lies in $\B(\C)$ is the relation $\sum
z_i\wedge(1-z_i)=0\in\C^*\wedge\C^*$ on the simplex parameters $z_i$.
For a more restrictive type of ideal triangulation than those
considered here this relation has been attributed to Thurston
(unpublished) by Gross \cite{gross} (according to \cite{zagier}). It
also follows easily from \cite{neumann-zagier} (see also
\cite{neumann}).  We give a cohomological proof here.

Recall (e.g., \cite{reid,neumann-reid1}) that the {\em invariant trace 
field} $k(M)$ is the field generated over $\Q$ by squares of
traces of elements of $\Gamma$.  It is known that modulo torsion the
Bloch group $\B(k)$ of a number field is isomorphic to $\Z^{r_2}$,
where $r_2$ is the number of complex embeddings of $k$, so
$\B(k)\otimes\Q\cong\Q^{r_2}$.  Moreover, if $k$ is given as a
subfield of $\C$, then naturally $\B(k)\otimes\Q\subset\B(\C)\otimes
\Q$.
\begin{theorem}\label{theorem2}
As an element of\/ $\B(\C)\otimes\Q$, the invariant\/ $\beta(M)$
lies in the subgroup\/
$\B(k(M))\otimes\Q$.
\end{theorem}
In the non-compact case we can define $\beta(M)$ directly as an
element of $\B(k(M))$ which is independent of triangulation, so the
``$\otimes\Q$'' of the previous theorem can be deleted, but we do not
know if it can in the compact case. Theorem \ref{restrict} is of
interest in this regard.

The Chern-Simons invariant $\CS(M)$ is determined modulo rational
multiples of $\pi^2$ by $\beta(M)$.  Chern and Simons defined what is
now called the Chern Simons invariant in (\cite{chern-simons}) for any
compact $(4n-1)$-dimensional Riemannian manifold. Meyerhoff
\cite{meyerhoff} extended the definition in the case of hyperbolic
3-manifolds to allow noncompact ones, that is hyperbolic 3-manifolds
with cusps. The Chern-Simons invariant $\CS(M)$ of such a hyperbolic
3-manifold $M$ is an element in $\R/\pi^2\Z$.  There is a map
$\rho\colon \B(\C)\to \C/\Q$ called the ``Bloch regulator map'' whose
definition we recall in section \ref{chern simons}.

\begin{theorem}\label{formula}
$\rho(\beta(M))=\frac i{2\pi^2}(\vol(M)+i\CS(M))\in\C/\Q$.
\end{theorem}

As is pointed out in \cite{neumann-yang}, Theorems \ref{theorem2} and
\ref{formula} put strong restrictions on $\CS(M)$. For example, it
follows that $\CS(M)$ is rational (by which we mean that it is zero in
$\R/\pi^2\Q$) if $k(M)$ is a quadratic extension of a totally real
field, and --- assuming the ``Ramakrishnan Conjecture'' --- $\CS(M)$ is
irrational if $k(M)\cap\overline{k(M)}\subset\R$.

We also discuss a definition of our invariant for any homomorphism
$f\colon \Gamma\to\PGL(2,\C)$.  It then generally lies in $\P(\C)$ rather
than $\B(\C)$, but it equals $\beta(M)$ for the discrete embedding of
$\Gamma$. We generalize this also to higher dimensions, but it is not
clear at this point what the significance of this is.

In section \ref{examples} we describe several interesting examples.
The final section \ref{appendix} is an appendix describing a scissors
congruence interpretation of the Bloch group $\B(\C)$.

Many of the results of this paper were announced in
\cite{neumann-yang-era}. The question of defining invariants in
$K$-groups for hyperbolic manifolds is also investigated by
A. Goncharov. We thank him for sharing his preprint \cite{goncharov}
with us.

\section{Preliminaries}
\subsection{Ideal simplices and degree one ideal triangulations}
\label{ideal triang}
We shall denote the standard compactification of $\H^3$ by $\overline
\H^3 = \H^3\cup\CP^1$. An ideal simplex $\Delta$ with vertices
$z_1,z_2,z_3,z_4\in\CP^1$ is determined up to congruence by the cross
ratio
$$z=[z_1:z_2:z_3:z_4]=\frac{(z_3-z_2)(z_4-z_1)}{(z_3-z_1)(z_4-z_2)}.$$
This $z$ lies in the upper half plane of $\C$ if the orientation
induced by the given ordering of the vertices agrees with the
orientation of $\H^3$. Permuting the vertices by an even (i.e.,
orientation preserving) permutation replaces $z$ by one of
$$
z,\quad 1-\frac 1z, \quad\text{or}\quad \frac 1{1-z},
$$
while an odd permutation replaces $z$ by
$$
\frac 1z, \quad\frac z{z-1},\quad\text{or}\quad 1-z.
$$
We will also allow degenerate ideal simplices where the vertices lie
in a plane, so the parameter $z$ is real.  However, we always require
that the vertices are distinct.  Thus the parameter $z$ of the simplex
lies in $\C-\{0,1\}$ and every such $z$ corresponds to an ideal
simplex.

If one takes finitely many geometric 3-simplices and glues them
together by identifying all the 2-faces in pairs then one obtains a
cellular complex $Y$ which is a manifold except possibly at isolated
points.  If the complement of these ``bad'' points is oriented we 
call $Y$ a {\em geometric 3-cycle}.  In this case the complement
$Y-Y^{(0)}$ of the vertices is an oriented manifold.

Suppose $M^3=\H^3/\Gamma$ is a hyperbolic manifold\footnote{Throughout
this paper ``hyperbolic manifold'' means oriented hyperbolic
3-manifold of finite volume.  Similarly Kleinian groups are assumed to
have finite covolume}.  A {\em degree one ideal triangulation
of $M$} consists of a geometric 3-cycle $Y$ plus a map $f\colon
Y-Y^{(0)}\to M$ satisfying
\begin{itemize} \item $f$ is degree one almost everywhere
in $M$; \item for each 3-simplex $S$ of $Y$ there is a map $f_S$ of
$S$ to an ideal simplex in $\overline\H^3$, bijective on vertices,
such that $f|S-S^{(0)}\colon S-S^{(0)}\to M$ is the composition
$\pi\circ (f_S|S-S^{(0)})$, where $\pi\colon \H^3\to M$ is the projection.
\end{itemize}

In \cite{thurston2} Thurston shows that any compact hyperbolic
3-manifold has degree one ideal triangulations with $|Y|\simeq
M$. Ideal triangulations also arise ``in practice'' (e.g., in the
program SNAPPEA for exploring hyperbolic manifolds --- \cite{weeks})
as follows.  Epstein and Penner in
\cite{epstein-penner} show that any non-compact $M$ has a genuine
ideal triangulation, that is, one for which $f$ is arbitrarily closely
deformable to a homeomorphism (they actually give an ideal polyhedral
subdivision; to subdivide these polyhedra into ideal tetrahedra 
it is conceivable that one
may need flat ideal tetrahedra to match triangulations of faces of
polyhedra --- see section \ref{appendix} for more details).  
The ideal simplices can be deformed to give degree one
ideal triangulations (based on the same geometric 3-cycle $Y$) on
almost all manifolds obtained by Dehn filling cusps of $M$ (see e.g.,
\cite{neumann-zagier}).

\subsection{Bloch group}\label{ssbloch} There are several different
definitions of the Bloch group in the literature. They differ at most
by torsion and they agree with each other for
algebraically closed fields.  We shall use the following.

\begin{definition}\label{def-bloch}
Let $k$ be a field.  The {\em pre-Bloch group $\P(k)$} is the quotient
of the free $\Z$-module $\Z (k-\{0,1\})$ by all instances of the
following relations:
\begin{gather}
[x]-[y]+[\frac yx]-[\frac{1-x^{-1}}{1-y^{-1}}]+[\frac{1-x}{1-y}]=0,
\label{5term}\\
[x]=[1-\frac 1x]=[\frac 1{1-x}]=-[\frac1x]=-[\frac{x-1}x]=-[1-x].
\label{invsim}
\end{gather}
The first of these relations is usually called the {\em five term
relation}.
The {\em Bloch group $\B(k)$} is the kernel of the map
$$
\P(k)\to k^*\wedge_\Z k^*,\quad [z]\mapsto 2(z\wedge(1-z)).
$$
\end{definition}

(The above five term relation is the one of Suslin \cite{suslin}. It
had misprints in its formulations in \cite{neumann-yang} and
\cite{neumann-yang-era}. Dupont and Sah use
$[x]-[y]+[y/x]-[(1-y)/(1-x)]+[(1-y^{-1})/(1-x^{-1})]=0$, which is
conjugate to Suslin's by the self-map $z\mapsto z^{-1}$ of
$\Z(k-\{0,1\})$.  One easily deduces that this map induces an
isomorphism between the Bloch groups resulting from the two choices of
five-term relation.  We will thus use Suslin's without further
comment.  The justification of the choice comes from Suslin's version
of the Bloch group which we do not discuss here, see \cite {suslin} or
\cite{neumann-yang} for more details.)

Dupont and Sah's definition of the Bloch group does not use the
relations (\ref{invsim}).  We shall need their versions later so we
will denote their groups obtained by omitting relations (\ref{invsim})
by $\P'(k)$ and $\B'(k)$ (this is the reverse of their convention).

Dupont and Sah show in \cite{dupont-sah} that $\P'(k)$ is more
natural than $\P(k)$ from a homological point of view. They also
show:
\begin{lemma}\label{p=p'} If
the characteristic of $k$ is not $2$ then $\P(k)$ and $\P'(k)$ differ
by at most torsion of order dividing $6$.  If, moreover, $k$ is
algebraically closed then $\P(k)=\P'(k)$.\qed
\end{lemma} 
Thus, if we are willing to ignore torsion or if we are working
over $\C$ we can use either definition.

\begin{definition}\label{def2}
Another definition starts
with a pre-Bloch group $\P(k)$ defined as the quotient of the free
$\Z-module$ $\Z (k\cup\{\infty\})$ by the following relations:
\begin{gather*}
[0]=[1]=[\infty]=0,\\
[x]-[y]+[\frac yx]-[\frac{1-x^{-1}}{1-y^{-1}}]+[\frac{1-x}{1-y}]=0.
\end{gather*}
It is not hard to see that this gives the same result as our definition.
\end{definition}
\begin{remark}\label{remark1}
For $k=\C$, the relations (\ref{invsim}) express the fact that
$\P(\C)$ may be thought of as being generated by isometry classes of
ideal hyperbolic 3-simplices.  The five term relation (\ref{5term})
then expresses the fact that in this group we can replace an ideal
simplex on four ideal points by the cone of its boundary to a fifth
ideal point.  As we show in an appendix (section \ref{appendix}), the
effect is that $\P(\C)=\P'(\C)$ is a group generated by ideal
polyhedra with ideal triangular faces modulo the relations generated
by cutting and pasting along such faces.
\end{remark}

\subsection{The Bloch invariant}\label{Bloch Invariant}

Let $f\colon Y-Y^{(0)}\to M$ be a degree 1 ideal triangulation of the
hyperbolic manifold $M$ as in subsection \ref{ideal triang} above.
Each 3-simplex of $Y$ maps to an ideal hyperbolic simplex in
$\overline\H^3$.  Let $z_1,\dots,z_n$ be the cross ratio parameters of
these ideal simplices.

\begin{definition} The {\em Bloch invariant $\beta(M)$} is the element
$\sum_1^n [z_j]\in\P(\C)$.  If the $z_j$'s all belong to a subfield
$K\subset \C$, we may consider $\beta(M)$ as an element of $\P(K)$.
\end{definition}

We show in section \ref{field sect} that if $\beta(M)$ can be defined
as above in $\P(K)$ then it lies in $\B(K)\subset\P(K)$ and is
independent of triangulation.

\section{Relative homology of $\Gamma$}\label{relative sect}

If $G$ is a group and $\Omega$ is a $G$-set then $\Z\Omega$ is a $\Z
G$-module.  Let $J\Omega$ be the kernel of the augmentation map
$\epsilon\colon \Z\Omega\to\Z$.  Then, following \cite{dicks-dunwoody}, we
define
$$H_n(G,\Omega)=H_n(G,\Omega;\Z)=\Tor_{n-1}^{\Z G}(J\Omega,\Z).$$
For our purposes it is convenient to capture the dimension shift in
this definition as follows.  Suppose
$$\cdots\to S_3\to S_2\to S_1\to J\Omega\to 0$$
is a $\Z G$-projective resolution of $J\Omega$.  Then
$H_n(G,\Omega)$ is the homology at index $n$ of the chain complex
$$\cdots\to (S_3)_G \to (S_2)_G
\to (S_1)_G\to 0,$$
where we are using the notation $$M_G:=M\otimes_{\Z G}\Z.$$

Let $S_n(\CP^1)$ denote the free abelian group generated by all
ordered 
$(n+1)$-tuples $\langle z_0,\ldots,z_n\rangle$ of {\em distinct}
points of $\CP^1$ modulo the relations
$$\langle z_0,\ldots,z_n\rangle=\sgn{\tau}\langle z_{\tau(0)},\ldots
,z_{\tau(n)}\rangle $$
for any permutation $\tau$ of $\{0,\ldots,n\}$.  With the standard
boundary map, they form a chain complex $S\subdot(\CP^1)$. This is
the cellular chain complex for the complete simplex on $\CP^1$
(considered as a discrete set), so it gives a resolution of $\Z$,
i.e., the following sequence
$$
\cdots\rightarrow S_2(\CP^1)\rightarrow S_1(\CP^1) \rightarrow
S_0(\CP^1)\rightarrow \Z\rightarrow 0
$$
is exact.

If $G\subset\PGL(2,\C)$ is any subgroup, then $G$ acts on $\CP^1$, so
the above sequence is a sequence of $\Z G$-modules.  We can truncate
this exact sequence to get an exact sequence
$$
\cdots\rightarrow S_2(\CP^1)\rightarrow S_1(\CP^1) \rightarrow
J\CP^1 \rightarrow 0.
$$
If this resolution of $J\CP^1$ were a free $\Z G$-resolution then
$S_{\bullet\ge1}(\CP^1)_G$ would compute the homology
$H_n(G,\CP^1)$.  Instead we only get a homomorphism
$$
H_n(G,\CP^1)\to H_n(S_{\bullet\ge1}(\CP^1)_G)
$$
for each $n$, determined by mapping any free resolution to the above
resolution. 

For the rest of the section, suppose we have a hyperbolic manifold
$M=\H^3/\Gamma$, so $\Gamma\subset\PGL(2,\C)$ is a torsion free
Kleinian group. Then $\Gamma$ acts on
$\CP^1=\bH[3]$. In this case we shall see that the
above homomorphism 
is an isomorphism in degree $3$ (and higher).

\begin{lemma}\label{free}
$S_n(\CP^1)$ is a free $\Z\Gamma$-module for $n\geq 2$.
\end{lemma}
\begin{proof} We will show that $\Gamma$ acts freely on the basis of
$S_n(\CP^1)$ for $n\geq 2$.  Since each element $\gamma \in \Gamma$ is
uniquely determined by its action on 3 distinct points of $\CP^1$, the
only way for $\gamma$ to fix a basis element of $S_n(\CP^1)$ for
$n\geq 2$ is if it acts as a permutation of $n+1$ distinct points in
$\CP^1$. But such a $\gamma$ must be a torsion element which is
impossible since $\Gamma$ is torsion free.
\end{proof}


\begin{proposition}\label{relative-isom}
{\rm 1.} $H_3(\Gamma, \CP^1)\cong\Z$.  
\item{\rm 2.} The above map $H_3(\Gamma, \CP^1)\to
H_3(S\subdot(\CP^1)_\Gamma)$ is an isomorphism.
\item{\rm 3.} Furthermore, if\/ $\Gamma$ is cocompact, then the natural map
$H_3(\Gamma;\Z)\to H_3(\Gamma,\CP^1)$ is also an isomorphism.
\end{proposition}
\begin{proof}
If $M$ is compact then part 1 follows from part 3, since
$H_3(\Gamma;\Z)=H_3(M;\Z)=\Z$.  Thus assume $M$ has cusps.
Let $\mathcal C$ be the set of {\em cusp points}, that
is, preimages of cusps of $M$ in $\bH[3] = \CP^1$, or, equivalently,
fixed points of parabolic elements of $\Gamma$.  Let $M_0$ be the
result of removing open horoball neighbourhoods of the cusps of $M$,
so $M_0$ is compact with toral boundary components.  Then $\calC$ can
also be identified with $\pi_0(\partial\widetilde{M_0})$, where
$\widetilde{M_0}$ is the universal cover. As described for instance in
\cite{dicks-dunwoody} (see also the proof of Lemma \ref{fund-class}),
one has $H_3(\Gamma, \mathcal C) \cong H_3(M_0,\partial M_0) \cong \Z$,
generated by the fundamental class $[M_0]$.  The short exact sequence
$$0\rightarrow J\calC \rightarrow J\CP^1\rightarrow \Z(\CP^1-\calC)
\rightarrow 0$$
gives rise to the long exact sequence
$$
\cdots \rightarrow H_3(\Gamma, \Z(\CP^1-\calC))\rightarrow H_3(\Gamma,\calC)
\rightarrow H_3(\Gamma, \CP^1) \rightarrow H_2(\Gamma,
\Z(\CP^1-\calC))\rightarrow \cdots.
$$
 By
Shapiro's Lemma (\cite{brown}) $H_i(\Gamma, \Z(\CP^1-\calC))$ is
isomorphic to the direct sum over the orbits of $\Gamma$ on
$\CP^1-\calC$ of $H_i$ of the isotropy groups of these orbits. Since
these isotropy groups are all trivial or infinite cyclic, both
$H_3(\Gamma, \Z(\CP^1-\calC))$ and $H_2(\Gamma, \Z(\CP^1-\calC))$ are
trivial. It follows that
$$H_3(\Gamma,\calC)\stackrel{\cong}\longrightarrow H_3(\Gamma,\CP^1),$$
completing the proof of part 1.

We prove parts 2 and 3 together.  Since $S\subdot(\CP^1)$ is a
resolution of $\Z$, the standard spectral sequence associated to
$H_\ast(\Gamma, S\subdot(\CP^1))$ converges to $H_\ast(\Gamma;
\Z)$. The $E^1_{p,q}$ term of the spectral sequence is $H_q(\Gamma,
S_p(\CP^1))$. It follows from Lemma \ref{free} and the fact that
$H_i(G;M)=\{0\}$ for any free $G$-module $M$ and $i\geq 1$, that
$E^1_{p,q} = \{0\}$ for $p\geq 2$. The isotropy group of the $\Gamma$
action of a point in $\CP^1$ can only be the trivial group, an
infinite cyclic group or a torus group (i.e., $\Z\oplus \Z$), which
all have homology dimension at most 2. By Shapiro's Lemma $H_i(\Gamma,
S_0(\CP^1))$ is therefore trivial for $i\geq 3$. The isotropy group in
$\Gamma$ of an unordered pair of distinct points of $\CP^1$ is either
trivial or an infinite cyclic group. Again by Shapiro's Lemma,
$H_i(\Gamma, S_1(\CP^1)) = 0$ for $i\geq 2$.

Putting everything together, we deduce that the $E^1$-terms of the
spectral sequence in total degree $\leq 4$ have the following picture
$$
\begin{matrix}
 0 \\
 0 & 0 \\
 H_2(\Gamma, S_0(\CP^1)) & 0 & 0 \\
 H_1(\Gamma, S_0(\CP^1)) & H_1(\Gamma, S_1(\CP^1)) & 0 & 0 \\
S_0(\CP^1)_\Gamma & S_1(\CP^1)_\Gamma &
S_2(\CP^1)_\Gamma & S_3(\CP^1)_\Gamma &
S_4(\CP^1)_\Gamma
\end{matrix}
$$

In fact, all the omitted terms are zero except those in the bottom
row. For the time being, assume the following

\begin{lemma}\label{d1}
The $d_1$ differential $ H_1(\Gamma, S_1(\CP^1)) \rightarrow
H_1(\Gamma, S_0(\CP^1)) $ is injective.
\end{lemma}
It then follows that the $d_2$ differential maps
$H_3(S\subdot(\CP^1)_\Gamma)$ to the trivial subgroup of
$H_1(\Gamma, S_1(\CP^1))$. Now if $\Gamma$ is cocompact, then it does
not contain any parabolic elements. Hence for any point in $\CP^1$ the
isotropy group in $\Gamma$ is either trivial or infinite cyclic. It
follows that $H_2(\Gamma, S_0(\CP^1))$ is trivial. So if $\Gamma$ is
cocompact, the spectral sequence collapses at $E^3$, and it follows
that
$H_3(\Gamma,\Z)\to H_3(S\subdot(\CP^1)_\Gamma)$
is an isomorphism.

To examine the relationship between
$H_3(S\subdot(\CP^1)_\Gamma)$ and
$H_3(\Gamma,\CP^1)$, we use the exact sequence
$$
\cdots\rightarrow S_2(\CP^1)\rightarrow S_1(\CP^1) \rightarrow
J\CP^1 \rightarrow 0
$$
The spectral sequence associated to $H_\ast(\Gamma, S_{\bullet\geq
1}(\CP^1))$ computes $H_\ast(\Gamma, \CP^1)$. The $E^1$-terms of this
spectral sequence are simply the $E^1$-terms of the spectral sequence
of $H_\ast(\Gamma, S\subdot(\CP^1))$ without the first column (with
the appropriate degree shift). The only possible non-zero $d_2$ is the
one from $H_3(S\subdot(\CP^1)_\Gamma)$ to $H_1(\Gamma,
S_1(\CP^1))$ which has already been shown to be trivial above. It
follows that this spectral sequence collapses at $E^2$ and therefore
$H_3(\Gamma, \CP^1)\to
H_3(S\subdot(\CP^1)_\Gamma)$ is an isomorphism.
\end{proof}

The proof of Lemma \ref{d1} requires a more geometric argument.
\begin{proof}[Proof of Lemma \ref{d1}]
We view $\CP^1$ as the natural boundary of $\H^3$.  Consider the
action of $\Gamma$ on the set $(\CP^1)^2-\Delta$ of ordered pairs
$(z_0,z_1)$ of distinct points of $\CP^1$.  If an element of $\Gamma$
takes one such pair $(z_0,z_1)$ to another $(z'_0,z'_1)$ then it takes
the $\H^3$-geodesic joining $z_0$ and $z_1$ to the $\H^3$-geodesic
joining $z'_0$ and $z'_1$.  Since $\Gamma$ has no torsion, no element
of $\Gamma$ can take $(z_0,z_1)$ to $(z_1,z_0)$. Let $X$ be a set of
orbit representatives of the $\Gamma\times C_2$ action on
$(\CP^1)^2-\Delta$, where $C_2$ is the order 2 group that maps
$(z_0,z_1)\mapsto (z_1,z_0)$. Then, as a $\Z\Gamma$-module,
$S_1(\CP^1)$ is the sum over $(z_0,z_1)\in X$ of $\Z\Gamma /
\Gamma_{(z_0,z_1)}$.  Let $X_0$ be the subset of $X$ consisting of
$(z_0,z_1)$ whose isotropy groups $\Gamma_{(z_0,z_1)}$ are non-trivial
(and hence infinite cyclic). Then by Shapiro's lemma
$$H_1(\Gamma, S_1(\CP^1))\cong\bigoplus\limits_{(z_0,z_1)\in X_0}
H_1(\Gamma,\Z(\Gamma/\Gamma_{(z_0,z_1)}))=\bigoplus\limits_{(z_0,z_1)\in
X_0} \Z=\Z X_0.
$$

Similarly, let $X_1$ denote the set of orbits in $\CP^1/\Gamma$ whose
isotropy groups are infinite cyclic. They are in fact the orbits of
points that appear in $X_0$. Again by Shapiro's Lemma, $H_1(\Gamma,
S_0(\CP^1))$ is $\Z X_1$. The differential $d_1$ takes $(z_0,z_1)\in\Z
X_0$ to $[z_1]-[z_0]$ in $\Z X_1$.

Now if $(z_0,z_1)\in X_0$ then $z_0$ and $z_1$ are in different
$\Gamma$-orbits, for if not then the $\H^3$-geodesic connecting $z_0$
and $z_1$ would map to a closed geodesic in $M=\H^3/\Gamma$ which is
asymptotic to itself with reversed direction, which is clearly absurd.
Similarly, if $(z_0,z_1)$ and $(z'_0,z'_1)$ are distinct elements of
$X_0$ then $z_0,z_1,z'_0,z'_1$ represent four distinct
$\Gamma$-orbits, for otherwise $(z_0,z_1)$ and $(z'_0,z'_1)$ would
represent two closed geodesics in $M$ which are asymptotic to each
other. The injectivity of $d_1$ now follows immediately.
\end{proof}

\section{The fundamental homology class}\label{fundamental sect}

Suppose $M=\H^3/\Gamma$ is a hyperbolic manifold with a degree one
ideal triangulation given by a geometric 3-cycle $Y$ and map
$f\colon Y-Y^{(0)}\to M$.  Form the pull-back covering
$$
\begin{CD}
\widehat{Y-Y^{(0)}} @>>> \H^3 \\
 @VVV             @VVV \\
Y-Y^{(0)} @>>> M
\end{CD}.
$$
We can complete to get a simplicial complex $\widehat Y$ with a
$\Gamma$-action (which is free except maybe at the vertices) plus a
$\Gamma$-equivariant map $\widehat Y\to \overline\H^3$.  Since this
map takes vertices to $\CP^1$ it induces a $\Gamma$-equivariant map of
$\widehat Y$ to the complete simplex on $\CP^1$.  Denote by
$C\subdot(\widehat Y)$ the simplicial chain complex
of $\widehat Y$. We get an induced map $C\subdot(\widehat Y)\to
S\subdot(\CP^1)$ of chain complexes and hence a map in homology
$$\gamma\colon H_\ast(Y)=H_\ast(C\subdot(\widehat Y)_\Gamma)\to
H_\ast(S\subdot(\CP^1)_\Gamma),$$
since $C\subdot(\widehat Y)_\Gamma$ is the simplicial chain complex
$C\subdot(Y)$ of $Y$.
There is also a natural homomorphism
$$
\mu\colon H_3(S\subdot(\CP^1)_\Gamma) \to \P(\C),
$$
given by sending any 4-tuple of points in $\CP^1$ to their
cross ratio. 
\begin{lemma}\label{fund-class0}
If $[Y]\in H_3(Y)$ is the fundamental class then
$\mu\circ\gamma[Y]\in\P(\C)$ is the Bloch invariant $\beta(M)$.
\end{lemma}
\begin{proof}
The fundamental class $[Y]$ is represented by the sum of the
3-simplices of $Y$.  Under $\mu\circ\gamma$ this maps to the sum of
the cross ratio parameters of the corresponding ideal simplices.  This
is $\beta(M)$ by definition.\end{proof}

Recall that in Proposition \ref{relative-isom} we showed a natural
isomorphism $H_3(\Gamma,\CP^1)\cong H_3(S\subdot(\CP^1)_\Gamma)$ and
showed both groups are infinite cyclic with a natural generator which
we will denote in both cases by $[M]$.
\begin{lemma}\label{fund-class}
$\gamma[Y]=[M]\in H_3(S\subdot(\CP^1)_\Gamma)$.
\end{lemma}
\begin{proof}
We introduce a homomorphism $v\colon H_3(S\subdot(\CP^1)_\Gamma)\to
\R$ and show first that $v(\gamma[Y])=\vol(M)$, the volume of $M$.  We
will then show the same for the generator of $H_3(\Gamma,\CP^1)$ to
complete the proof.

We will use $\langle z_0,\ldots,z_n\rangle _\Gamma$ to denote the
image of $\langle z_0,\ldots,z_n\rangle $ in $S_n(\CP^1)_\Gamma$. We
can think of $\langle z_0,\ldots,z_3\rangle _\Gamma$ as representing
an ideal simplex in $\H^3$ that is well defined up to the action of
$\Gamma$.  We define $v(\langle z_0,\ldots,z_3\rangle)$ to be plus or
minus the hyperbolic volume of this ideal simplex, with the sign
chosen according as the orientation of the simplex agrees or not with
the orientation of $\H^3$ (if the simplex is planar then the volume is
zero and orientation is irrelevant).  Given a 3-cycle $\alpha = \sum_i
n_i \langle z^i_0,z^i_1,z^i_2,z^i_3\rangle _\Gamma$ in
$S_3(\CP^1)_\Gamma$ we define $v(\alpha)= \sum_i n_i v( \langle
z^i_0,z^i_1,z^i_2,z^i_3\rangle)$. If we have five distinct points
$z_0,\ldots,z_4\in\CP^1$ then it is easy to see geometrically that
$v(\partial\langle z_0,\ldots,z_4\rangle)=0$.  It follows that $v$
induces a map, which we also call $v$:
$$
v\colon H_3(S\subdot(\CP^1)_\Gamma) \rightarrow \R.
$$
The value of $v$ on a class given by a degree one ideal triangulation
is the sum of the signed volumes of the ideal tetrahedra of the
triangulation, which is just the volume of $M$, by the degree one
condition and Fubini's theorem.

Now let $Z$ be the end compactification of $M$.  Then $Z=X/\Gamma$
where $X=\H^3\cup\calC$ with $\calC$ as in the previous section.  Note
that $Z$ is homeomorphic to the result of collapsing each boundary
component of $M_0$ to a point, where $M_0$ is as in the proof of
Proposition \ref{relative-isom}.  In particular,
$H_3(Z)=H_3(M_0,\partial M_0)=\Z$.

Consider a triangulation of $Z$ and the lifted triangulation of
$X$. Let $C\subdot(X)$ denote the simplicial chain complex of $X$.  We
can think of $S_q(\CP^1)$ as being generated by arbitrary
$(q+1)$-tuples of points of $\CP^1$ modulo the relations
$$\langle z_0,\ldots,z_q\rangle=\sgn(\tau)\langle z_{\tau(0)},\ldots
,z_{\tau(q)}\rangle $$
for any permutation $\tau$ of $\{0,\ldots,q\}$ and
$$\langle z_0,\ldots,z_q\rangle=0\quad\hbox{if the $z_i$ are not
distinct}.$$
We can thus map $C\subdot(X) \to S\subdot(\CP^1)$ by taking any
equivariant map of the vertices of the triangulation of $X$ to
$\CP^1$.  We claim that the induced map
$C\subdot(X)_\Gamma=C\subdot(Z)\to S\subdot(\CP^1)_\Gamma$ induces the
isomorphism $H_3(\Gamma,\calC)=H_3(M_0,\partial M_0)=H_3(Z)\to
H_3(S\subdot(\CP^1)_\Gamma)$ of Proposition \ref{relative-isom}.

If $M$ is compact then, since $X=\H^3$ is contractible, $C\subdot(X)$
is a free $\Z\Gamma$ resolution of $\Z$, and the above map is indeed
the map of Proposition \ref{relative-isom}.

Suppose $M$ is non-compact, so $\calC$ is non-empty. Let $\widetilde
C\subdot(X)$ be the reduced chain complex (so $\widetilde
C_q(X)=C_q(X)$ for $q>0$ and $\widetilde C_0(X)$ is the kernel of
augmentation $\epsilon\colon C_0(X)\to \Z$).  Then, since $X$ is
contractible,
$$\cdots \to C_2(X)\to C_1(X)\to \widetilde C_0(X)\to 0$$
is exact.  It is clearly a free $\Z\Gamma$-resolution of $\widetilde
C_0(X)$.  Moreover, $\widetilde C_0(X)$ is isomorphic to $J\calC\oplus
F$ in the notation of the previous section, where $F$ is the submodule
of $C_0(X)$ generated by {\em finite} vertices. Since $F$ is clearly a
free $\Z\Gamma$-module, it follows that $C_{\bullet\ge1}(X)_\Gamma=
C_{\bullet\ge1}(Z)$ computes $H_\ast(\Gamma,\calC)$. This gives the
isomorphism $H_3(M_0,\partial M_0)=H_3(Z)\cong H_3(\Gamma,\calC)$ used
in the proof of Proposition \ref{relative-isom}.

Now consider the above map $C\subdot(X)\to S\subdot(\CP^1)$.  The
induced map of reduced groups in degree zero, $\widetilde C_0(X) \to
J\CP^1$, is, up to a map of free summands, the inclusion $J\calC\to
J\CP^1$.  The map $C_{\bullet\ge1}(X)_\Gamma \to
S_{\bullet\ge1}(\CP^1)_\Gamma$ thus induces the map
$H_3(\Gamma,\calC)= H_3(Z) \to H_3(S\subdot(\CP^1)_\Gamma)$ of
Proposition \ref{relative-isom} as claimed.

To compute $v$ of the generator of $H_3(S\subdot(\CP^1)_\Gamma)$ we
must thus map the vertices of $X$ equivariantly to $\CP^1$ and then
sum the volumes of the ideal simplices in $\H^3$ corresponding to a
set of $\Gamma$-orbit representatives of the 3-simplices of $X$.
Replace the simplices of the triangulation of $X$ by geodesic
simplices.  Then the sum of the signed volumes of a set of orbit
representatives of these simplices is $\vol(M)$.  In fact, if we
consider the triangulation of $Z$ to be given by a map from an
abstract simplicial complex to $Z$, then the above sum of signed
simplex volumes is just the integral of the pull-back of the volume
form on $Z$ to this simplicial complex.  If we now move the vertices
in $\H^3$ of the triangulation of $X$ continuously and equivariantly
within $\H^3$, then we are just homotoping the map of the simplicial
complex to $Z$, so the sum of signed simplex volumes stays equal to
$\vol(M)$.  If we now let the vertices move continuously all the way
out to $\CP^1$, then the simplex volumes change continuously, so their
sum is still $\vol(M)$ when the vertices reach $\CP^1$.
\end{proof}

Combining the isomorphism of Proposition \ref{relative-isom} with the
above map $\mu$ gives a map
$$H_3(\Gamma,\CP^1)\to \P(\C).$$
The preceding lemmas imply immediately:
\begin{proposition}\label{fund-beta}
This homomorphism maps $[M]\in H_3(\Gamma,\CP^1)$ to $\beta(M)$.  In
particular $\beta(M)$ is independent of triangulation.\qed
\end{proposition}

\begin{remark}\label{maptoPC}
Note that the above map $H_3(\Gamma,\CP^1)\to \P(\C)$ is defined for
{\em any} subgroup $\Gamma\subset\PGL(2,\C)$, including $\PGL(2,\C)$
itself.  Now $\PGL(2,\C)$ acts transitively on the basis of
$S_2(\CP^1)$, so the boundary map $S_3(\CP^1)_{ \PGL(2,\C)}\to
S_2(\CP^1)_{\PGL(2,\C)}$ is trivial.  It follows that
$$H_3(S\subdot(\CP^1)_{\PGL(2,\C)})\stackrel{\cong}\longrightarrow\P(\C).$$
The map of the above proposition can thus be described as the map
$H_3(S\subdot(\CP^1)_\Gamma)\to H_3(S\subdot(\CP^1)_{\PGL(2,\C)}) =
\P(\C)$ induced by the inclusion $\Gamma\to \PGL(2,\C)$.
\end{remark}

\section{Completion of Proof of Theorem \ref{theorem1}}\label{proof sect}

For Theorem \ref{theorem1} it remains to show that $\beta(M)$ is in
$\B(\C)\subset\P(\C)$.  

We first assume $M$ is not compact.  The inclusion $J\calC\to\Z\calC$
induces (because of the dimension shift in the definition of
$H_\ast(\Gamma,\calC)$) a map $H_3(\Gamma,\calC)\to
H_2(\Gamma,\Z\calC)$.  Similarly, the inclusion $J\CP^1\to \Z\CP^1$
induces a map $H_3(\PGL(2,\C),\CP^1)\to H_2(\PGL(2,\C),\Z\CP^1)$.  We
thus have the following commutative diagram
$$\begin{CD}
H_3(\Gamma,\calC) @>>> H_2(\Gamma, \Z\calC)\\
@VVV  @VV\phi V \\
H_3(\PGL(2,\C),\CP^1) @>>> H_2(\PGL(2,\C), \Z\CP^1)
\end{CD}$$

Since $\PGL(2,\C)$ acts transitively on $\CP^1$ with isotropy group
the Borel subgroup $B$ of upper triangular matrices, we have
$H_2(\PGL(2,\C),\Z\CP^1)\cong H_2(B;\Z)$. Let $T$ denote the maximal
torus in $\PGL(2,\C)$ and $U$ the upper unipotent matrices.  Then
$H_2(B,\Z)\cong H_2(T,\Z)\cong \wedge^2\C^\ast$, induced by
$B\rightarrow B/U\cong T$, (cf. \cite{dupont-sah} (A11)).  By
\cite{dupont-sah}, one has the following commutative diagram
$$
\begin{CD}
H_3(\PGL(2,\C),\CP^1) @>>> H_2(\PGL(2,\C), \Z\CP^1)\\
@VV\cong V @VV\cong V \\
\P(\C) @>\lambda >> \wedge^2\C^\ast
\end{CD}
$$
and the map $\lambda$ is simply given by $\lambda[z] = 2 (z\wedge
(1-z))$. 

So to show that the element $\beta(M)$ is in fact in $\B(\C)$, it
suffices to show that the homomorphism $\phi$ is trivial.  By
Shapiro's lemma, $H_2(\Gamma,\Z\calC)$ is the sum of $H_2$ of isotropy
groups of $\Gamma$ on $\calC$. These isotropy groups are unipotent
subgroups of $\PGL(2,\C)$, so $H_2$ of such an isotropy group when
mapped to $H_2(B,\Z)$ under $\phi$ factors through $H_2(U,\Z)$. Since
$H_2(B,\Z)\cong H_2(T,\Z)$ is induced by $B\rightarrow B/U\cong T$,
 the homomorphism from $H_2(U,\Z)$ to
$H_2(B,\Z)$ is trivial. It follows that $\phi$ is trivial, hence
$\beta[M]$ is an element in $\B(\C)$. 

If $M$ is compact then an argument was given in the Introduction.  We
expand on it since we need it again below. In this case $\beta(M)$ is
in the image of $H_3(\Gamma)\to H_3(\PGL(2,\C),\CP^1)\to\P(\C)$,
which factors through $H_3(\PGL(2,\C))$.  Since the sequence
$H_3(\PGL(2,\C))\to H_3(\PGL(2,\C),\CP^1)\to H_2(\PGL(2,\C),\Z\CP^1)$
is exact, the result again follows from the above diagram. This
finishes the proof of Theorem \ref{theorem1}.\qed

\section{Restricting the field of definition}\label{field sect} 

There are two situations in which we can define the Bloch invariant
$\beta(M)$ as an element of $\P(K)$ for a subfield $K\subset \C$
rather than as an element of $\P(\C)$.

\smallbreak\noindent{\bf Case 1.}
If $M$ is a hyperbolic manifold with a degree one ideal triangulation
into ideal simplices $\Delta_1,\dots,\Delta_n$ we call the subfield of
$\C$ generated by the cross ratio parameters $z_1,\dots,z_n$ of these
simplices the {\em tetrahedron field} associated with the
triangulation.  If this tetrahedron field is a subfield of $K$ we can
define $$\beta(M) :=\sum_{i=1}^n[z_i]\in\P(K).$$

\smallbreak\noindent{\bf Case 2.}
If $M=\H^3/\Gamma$ with $\Gamma\subset\PGL(2,K)$ then we can define
$\beta(M)$ as the image of $[M]$ under the composition
$$H_3(\Gamma,\calC)\stackrel\cong\to H_3(S\subdot(\CP^1(K))_\Gamma)\to
H_3(S\subdot(\CP^1(K))_{\PGL(2,K)})\stackrel{\cong}\to\P(K),$$ where
$\CP^1(K)=K\mathbb P^1$ is the set of $K$-rational points of $\CP^1$.
\smallbreak

\begin{theorem}\label{restrict}
If $M=\H^3/\Gamma$ has a tetrahedron field contained in $K$, then
$\Gamma$ has a discrete embedding into $\PGL(2,\C)$ with image in
$\PGL(2,K)$.  That is, if Case 1 holds then so does Case 2. Moreover,
in either case\/ $\beta(M)$ lies in\/ $\B(K)$ and only depends on $M$
and the field $K$.
\end{theorem}

\begin{proof}
If $f\colon Y-Y^{(0)}\to M$ is a degree one ideal triangulation then,
as above, we can lift to a map $\widehat{Y}\to \overline\H^3$.  Let
$V\subset\CP^1=\bH[3]$ be the image of the set of vertices of
$\widehat Y$.  Then, as in \cite{neumann-reid1}, if we apply an
isometry to put three of the points of $V$ at $0$, $1$, and $\infty$,
then the fact that $Y$ is connected implies that $V$ will be a
subset of $\CP^1(K)$, where $K$ is the field generated by the cross
ratio parameters of the simplices of $K$.  Since $V$ is a
$\Gamma$-invariant set, it also follows as in \cite{neumann-reid1}
that $\Gamma$ is then in $\PGL(2,K)$.  We can now repeat the proof of
Theorem \ref{theorem1} using $\CP^1(K)$ in place of $\CP^1$.  The
proof that $\beta(M)\in\P(K)$ is independent of choices is just as
before.  The proof that it lies in $\B(K)$ needs slightly more care.

We first assume $M$ is non-compact. We consider the versions of the
commutative diagrams of section \ref{proof sect} with $\C$ replaced by
$K$.  The first,
$$\begin{CD}
H_3(\Gamma,\calC) @>>> H_2(\Gamma, \Z\calC)\\
@VVV  @VV\phi V \\
H_3(\PGL_2(K),\CP^1(K)) @>>> H_2(\PGL_2(K), \Z\CP^1(K)),
\end{CD}$$
is unproblematic.  For the second, one must follow the argument
presented in the appendix of Dupont-Sah \cite{dupont-sah} carefully.
It gives a commutative diagram
$$
\begin{CD}
H_3(\PGL_2(K),\CP^1(K)) @>>> H_2(\PGL_2(K), \Z\CP^1(K))\\
@VV{T}V @VV\cong V \\
\P'(K) @>\lambda >> \wedge^2K^\ast
\end{CD}
$$ with the only difference to the diagram for $K=\C$ being that $T$
is now only an isomorphism modulo 2-torsion and we are using the
Dupont-Sah version $\P'(K)$ of the Bloch group rather than $\P(K)$
(see section \ref{ssbloch}). The map $\lambda$ is still $[z]\mapsto
2(z\wedge(1-z))$. Since this vanishes on the extra relations
(\ref{invsim}) that define $\P(K)$ from $\P'(K)$, we can replace
$\P'(K)$ by $\P(K)$ in the above diagram, and the rest of the argument
carries through as in the case $K=\C$.  Similarly, the argument in the
compact case also carries through.
\end{proof}

Now if $M$ has cusps, then by \cite{epstein-penner}, $M$ has a genuine
decomposition into convex ideal polyhedra with vertices at the cusps.
We can further subdivide these polyhedra into ideal tetrahedra. These
subdivisions may not agree on common faces of the ideal polyhedra, in
which case we must add some flat ideal tetrahedra to mediate between
the two triangulations of the faces in question.  We obtain what we
call a ``genuine'' ideal triangulation.  In \cite{neumann-reid1}, it
was shown that all the cross ratio parameters of these ideal
tetrahedra lie in the invariant trace field $k=k(M)$. Hence
$\beta(M)$ can in fact be defined in $\P(k)$.

\begin{corollary} If $M$ is non-compact then $\beta(M)$ is
well-defined in $\B(k)$, where $k=k(M)$ is the invariant trace field
of $M$.\qed\end{corollary}

\begin{proof}[Proof of Theorem \ref{theorem2}] 
If $\Gamma$ can be conjugated into $\PGL(2,K)$, then $K$ is called a
potential coefficient field of $\Gamma$ (or $M$) in
\cite{macbeath}. By theorem \ref{restrict}, $\beta(M)$ can be defined in
$\B(K)$ for each potential coefficient field $K$. According to
\cite{macbeath}, one can find infinitely many such potential
coefficient fields $K$, which are at most quadratic extensions of the
trace field $K_0$ of $\Gamma$. Choose two of them, say $K_1$ and
$K_2$.  Then $K_1\cap K_2=K_0$. Let $L$ be any field which contains
$K_1$ and $K_2$.  The inclusions $K_i\to L$ induce modulo torsion
injections $\B(K_i)\to\B(L)$ for $i=0,1,2$ (see e.g.,
\cite{neumann-yang}). Since we are willing to work modulo torsion, we
shall identify each $\B(K_i)$ with its image in $\B(L)$.  In
\cite{neumann-yang} it is shown that with these identifications
$\B(K_0)\subset\B(K_1)\cap \B(K_2)$ with torsion quotient. By Theorem
\ref{restrict} we know
$\beta(M)\in\B(L)$ is in the subgroup $\B(K_1)\cap\B(K_2)$ and hence
some positive multiple is in the subgroup $\B(K_0)\subset\B(L)$. In
particular, $\beta(M)$ is well defined in
$\B(\tracefield(M))\otimes\Q$.

Now, if $\Gamma^{(2)}$ is the subgroup of $\Gamma$ generated by
squares of elements in $\Gamma$ then $\Gamma/\Gamma^{(2)}$ is an
elementary abelian group of order $2^s$ for some $s$.  In \cite{reid}
it is shown that the trace field of $\Gamma^{(2)}$ is the invariant
trace field $k(M)$ of $M$.  Let $M^{(2)}=\H^3/\Gamma^{(2)}$.  Then
clearly $\beta(M^{(2)})=2^s\beta(M)$ and since $\beta(M^{(2)})\in
\B(k(M))$ it follows that there is a well defined element
$\beta_{k(M)}(M) \in \B(k(M))\otimes\Q$ whose image is $\beta(M)$ in
$\B(\C)\otimes\Q$.
\end{proof}

\section{Chern-Simons invariant}\label{chern simons}

Theorem \ref{formula} involved the {\em Bloch
regulator map}
$$\rho\colon \B(\C) \longrightarrow \C/\Q,$$
which is defined as follows. For $z\in \C-\{0,1\}$, define
$$
\rho(z) = \frac{\log z}{2\pi i}\wedge\frac{\log (1-z)}{2\pi i} +
1\wedge \frac{{\mathcal R}(z)}{2\pi^2},
$$
where ${\mathcal R}(z)$ is the ``Rogers dilogarithm function''
$${\mathcal R}(z)=\frac12\log(z)\log(1-z)-\int_0^z\frac{\log(1-t)}tdt.$$
See section 4 of \cite{dupont-sah} or \cite{hain} for details on how
to interpret this formula.  This $\rho$ vanishes on the relations
(\ref{5term}) and (\ref{invsim}) which define $\P(\C)$ and hence
$\rho$ induces a map
$$\rho\colon \P(\C)\longrightarrow\C\wedge_\Z\C.$$
This fits in a commutative diagram
$$\begin{matrix}
\P(\C)&\stackrel{\mu}{\longrightarrow}&\C^*\wedge\C^*\\
\downarrow\scriptstyle\rho&&\downarrow\scriptstyle=\\
\C\wedge\C&\stackrel{\epsilon}{\longrightarrow}&\C^*\wedge\C^*\\
\end{matrix}
$$ where $\epsilon=2(e\wedge e)$ with $e(z)=\exp(2\pi iz)$. 
The kernel of $\mu$ is $\B(\C)$ and the kernel of $\epsilon$ is
$\C/\Q$. Hence $\rho$ restricts to give the desired map
$\rho\colon\B(\C)\to\C/\Q$.

Recall that Theorem \ref{formula} is the formula
$$
\frac{2\pi^2}i\rho(\beta(M))=\vol(M)+i\CS(M)\in\C/(i\pi^2\Q).
$$
The volume part of this result is not hard (see e.g.
\cite{dupont-sah}, \cite{neumann}, \cite{neumann-yang}).  In fact the
imaginary part of $2\pi^2\rho$ can be extended to a map
$\vol\colon\P(\C)\to \R$ given on generators by $[z]\mapsto D_2(z)$
where $D_2$ is the {\em Bloch-Wigner
dilogarithm} (cf. \cite{bloch})
$$D_2(z) = \Im \ln_2(z) + \log |z|\arg(1-z),\quad z\in \C -\{0,1\}.$$
The name $\vol$ is justified because $D_2(z)$ is the hyperbolic volume
of an ideal tetrahedron $\Delta$ with cross ratio $z$. 

The real part of the formula, giving Chern-Simons invariant, lies
deeper.  If $M=\H^3/\Gamma$ is compact then Dupont \cite{dupont1}
proved\footnote{The proof in \cite{dupont1} is for the ``flat''
Chern-Simons invariant and the assumption implicit there that this
agrees with the ``riemmanian'' Chern-Simons invariant for $M$ is
confirmed in \cite{dupont-kamber}.} the above
formula with $\beta(M)$ replaced by the image of the fundamental class
$[M]\in H_3(\Gamma;\Z)$ under $H_3(\Gamma;\Z)\to
H_3(\PGL(2,\C)^\delta;\Z)\to \B(\C)$.  But this image is $\beta(M)$ by
Propositions \ref{relative-isom} and \ref{fund-beta}.  Thus Theorem
\ref{formula} follows in this case.  If $M$ is non-compact it follows
from the formula for Chern-Simons invariant of \cite{neumann}.  We
give some of the details for completeness.

Suppose $M$ has an ideal triangulation which subdivides it into $n$
ideal tetrahedra
$$
M = \Delta_1 \cup\dots\cup \Delta_n.
$$
Choose an ordering of the vertices of the $j$-th tetrahedron that is
compatible with the orientation of $M$ and let $z_j^0$ be the
cross ratio 
parameter which then describes this tetrahedron. Let
$$
\calZ^0=\left(\begin{matrix}
\log z^0_1\\
\vdots\\
\log z^0_n\\
\log(1-z^0_1)\\
\vdots\\
\log(1-z^0_n)\end{matrix}\right).
$$
Recall from \cite{neumann-zagier} (see also \cite{neumann} which we
are following here) 
that if $M$ has $h$ cusps then the $z^0_j$ are
determined by so-called consistency and cusp relations which can be
written in the form 
$$
U\calZ^0 = \pi i{\bf d},
$$
where $U$ is a certain integral $(n+2h)\times2n$-matrix and 
$$
{\bf d}=\left(\begin{matrix}
d_1\\
\vdots\\
d_{n+2h}\end{matrix}\right)
$$
is some integral vector.  Geometrically, the consistency relations say
that the tetrahedra fit together around each edge of the triangulation
and the cusp relations say that generators of the cusp groups
represent parabolic isometries.

The consistency relations are given by a $n\times 2n$ submatrix of $U$
and the component containing $z^0=(z^0_1,\dots,z^0_n)$ in the set of
$z=(z_1,\dots,z_n)$ that satisfy these relations will be called {\em
Dehn surgery space} and denoted $\mathcal D$ (this is actually a
$2^h$-fold branched cover of what is usually called Dehn surgery
space, see \cite{neumann-zagier}, but the difference is irrelevant to
the current discussion).

The equation
\begin{equation}\label{Uc=d}
U{\bf c}={\bf d}
\end{equation}
has a solution ${\bf c}=\calZ^0/\pi i \in \C^{2n}$.
Since $U$ is an integral matrix, equation (\ref{Uc=d}) also has solutions
$$
{\bf c}=\left(\begin{matrix}
c_1'\\
\vdots\\
c_n'\\
c_1''\\
\vdots\\
c_n''
\end{matrix}\right)\in \Q^{2n}.
$$
In \cite{neumann} it is shown that solutions ${\bf c}$ can be found in
$\Z^{2n}$, in fact even in a certain affine sublattice of $\Z^{2n}$.

Let $M'$ be the result of a hyperbolic Dehn filling on $M$ obtained by
deforming the parameter $z^0=(z^0_1,\dots,z^0_n)$ to a new value
$z=(z_1,\dots,z_n)$ in Dehn surgery space $\mathcal D$ (cf.~e.g.,
\cite{neumann-zagier}).  Topologically $M'$ differs from $M$ in that a
new closed geodesic $\gamma_j$ has been added at the $j$-th cusp for
some $j\in\{1,\dots,h\}$.  Let $\lambda_j$ be the complex number which
has real part equal to the length of this geodesic and imaginary part
equal to its torsion (the latter is only well-defined modulo $2\pi$).
If no geodesic has been added at the $j$-th cusp we put $\lambda_j =
0$.

\begin{theorem}[\cite{neumann}]\label{cs-formula}
Given any solution ${\bf c}\in\Q^{2n}$ to equation\/ {\rm(\ref{Uc=d})}, there
exists a constant $\alpha=\alpha({\bf c})\in i\pi^2\Q$ such that if $M'$
is any result of hyperbolic Dehn filling obtained by deforming $M$ as
above, then
$$
\vol(M')+i\CS(M')=\alpha-\frac{\pi}2\sum^h_{j=1}\lambda_j-i\sum^n_{\nu=1}
\biggl({\mathcal R}(z_\nu)-\frac{i\pi}2\bigl(
{c_\nu'}\log(1-z_\nu)-{c_\nu''}\log(z_\nu)\bigr)\biggr)
.
$$
Moreover, if ${\bf c}\in\Z^{2n}$ then conjecturally
$\alpha\in\frac{i\pi^2}{12}\Z$, while if ${\bf c}$ is in the sublattice
mentioned above then conjecturally
$\alpha\in\frac{i\pi^2}{6}\Z$. 
\end{theorem}

The following summary of the proof in \cite{neumann} also explains how
Theorem \ref{formula} follows.  

In \cite{yoshida} Yoshida proved a conjecture of \cite{neumann-zagier}
that a formula of the following form hold on the space of Dehn
fillings of $M$:
$$\vol(M')+i\CS(M')=-\frac{\pi}2\sum^h_{j=1}\lambda_j +f(z),$$
where $f(z)$ is some analytic function of $z=(z_1,\dots,z_n)\in\mathcal
D$.  The formula of the above theorem has this form and it is not hard
to verify that its real part is correct with $\alpha=0$.  We thus have
two analytic functions whose real parts agree at the points $z\in\mathcal
D$ that correspond to Dehn fillings.  These points limit on $z^0$,
which is a smooth point of $\mathcal D$, from all tangent directions. It
follows that the two analytic functions agree up to an imaginary
constant on $\mathcal D$.  This imaginary constant is $\alpha$.  On the
other hand, it is shown in \cite{neumann} that the right side of the formula 
of Theorem \ref{cs-formula} without the constant $\alpha$ gives
$\frac{2\pi^2}i \rho(\beta(M'))$ 
modulo $i\pi^2\Q$.  We have already shown this equals
$\vol(M')+i\CS(M')$ modulo $i\pi^2\Q$ if $M'$ is compact.  Thus 
$\alpha\in$ $i\pi^2\Q$, so Theorem \ref{cs-formula} is 
proved.  (Note that we have used that Theorem
\ref{formula} holds in the compact case, which Dupont \cite{dupont1}
proved for the image
of the fundamental class rather than for $\beta(M)$.  The result 
that the image of fundamental class equals 
$\beta(M)$ (Proposition \ref{fund-beta})
thus fills a gap in the proof of Theorem \ref{cs-formula}
in \cite{neumann}.)

We have just mentioned that \cite{neumann} shows that the right side
of the formula of Theorem \ref{cs-formula} equals
$\frac{2\pi^2}i\rho(\beta(M)$ modulo $i\pi^2\Q$. Applying this and
Theorem \ref{cs-formula} to $M$ itself proves Theorem \ref{formula} in
the non-compact case.

\section{Generalizations: Higher Dimensions and Homomorphisms}

\begin{definition}
Let $S_q(\bH)$ be the abelian group generated by arbitrary
$(q+1)$-tuples of points of $\bH$ modulo the relations
$$\langle z_0,\ldots,z_q\rangle=\sgn{\tau}\langle z_{\tau(0)},\ldots
,z_{\tau(q)}\rangle $$
for any permutation $\tau$ of $\{0,\ldots,q\}$ and
$$\langle z_0,\ldots,z_q\rangle=0\quad\hbox{if the $z_i$ are not
distinct}.$$  Define a boundary map $S_q(\bH)\to S_{q-1}(\bH)$ by
the usual formula $\partial\langle z_0,\dots,z_q\rangle=\sum_{i=0}^q
(-1)^i \langle z_0,\dots,\widehat{z_i},\dots,z_q\rangle$.  
Note that $S\subdot(\bH)_{\Isom^+(\H^n)}$ is the result of adding
the relations
$$\langle gz_0,\dots,gz_q\rangle=\langle z_0,\dots,z_q\rangle$$ for
$g\in\Isom^+(\H^n)$ to the above definition. We define
$$\P_n:=H_n(S\subdot(\bH)_{\Isom^+(\H^n)}).$$ 
In particular, if $n=3$ then $S\subdot(\bH[3])=S\subdot(\CP^1)$ and
$\P_3=\P(\C)$.
\end{definition}

Now let $M^n$ be a manifold which is homeomorphic to the interior of a
compact manifold $M_0$ with (possibly empty) boundary such that the
universal cover $\widetilde{M_0}$ and all its boundary components are
contractible.  For example, a complete hyperbolic $n$-manifold of
finite volume has this property.  Let $\Gamma=\pi_1(M)=\pi_1(M_0)$. We
will define an invariant of a homomorphism $f\colon\Gamma\to
\Isom^+(\H^n)$ which generalizes the invariant $\beta(M)$ of previous
sections.  We shall need the homomorphism to satisfy a condition which
we describe and discuss below.

$\Gamma$ acts by covering transformations on $\widetilde{M_0}$.  Let
$X$ be the end compactification of $\widetilde M$ and $Z$ the end
compactification of $M$ (these can be obtained by collapsing each
boundary component of $\widetilde{M_0}$ respectively $M_0$ to a
point).  Denote $\calC=X-\widetilde M$.  For each $c\in \calC$ denote
$P_c=\{g\in \Gamma:gc=c\}$.  $P_c$ is isomorphic to the fundamental
group of the boundary component of $M_0$ corresponding to $c$.

\begin{condition}\label{condition}
We assume that $f\colon\Gamma\to \Isom^+(\H^n)$ has the property that
$f(P_c)$ fixes some point $x_c\in\bH$ for each $c$.  As we discuss at
the end of this section, this condition can be relaxed and is then
automatically satisfied in many cases, for instance if $M$ is an
odd-dimensional hyperbolic manifold.
\end{condition}

$\Gamma$ acts on $\bH$ via the homomorphism $f$. We can choose the
assignment $h\colon c\mapsto x_c$ to be $\Gamma$-equivariant, since if
$f(P_c)$ fixes $x_c$ then $f(P_{gc})=f(gP_cg^{-1})$ fixes $f(g)x_c$.
We then get an induced map $H_n(\Gamma,\calC)\to
H_n(\Isom^+(\H^n),\bH)$ and the image $\beta_h(f)\in
H_n(\Isom^+(\H^n),\bH)$ of the fundamental class in
$H_n(\Gamma,\calC)=H_n(M_0,\partial M_0)=\Z$ is an invariant of the
given situation. If $M$ is compact we use instead the image of the
fundamental class under $H_n(\Gamma)\to H_n(\Isom^+(\H^n))\to
H_n(\Isom^+(\H^n),\bH)$. 

In general $\beta_h(f)$ presumably depends on the choice of $h$. 
But trivially:

\begin{proposition}
If $h$ is unique, for example if $M$ is compact, or $M$ is a finite
volume hyperbolic manifold and $f\colon\Gamma\to \Isom^+(\H^n)$ the
homomorphism that determines its hyperbolic structure, then\/
$\beta_h(f)\in H_n(\Isom^+(\H^n),\bH)$ 
is a well defined invariant of $f$ that generalizes the\/
$\beta(M)$ of previous sections.
\qed\end{proposition}

We also have a natural map
$$\mu\colon H_n(\Isom^+(\H^n),\bH)\to \P_n$$ 
generalizing the map
$H_3(\PGL(2,\C),\CP^1)\stackrel\cong\longrightarrow \P(\C)$ of remark
\ref{maptoPC}.

\begin{theorem}\label{beta of homom} 
$\mu(\beta_h(f))$ does not depend on $h$ and is thus
an invariant just of $M$ and $f\colon\pi_1(M)\to\Isom^+(\H^n)$.  
We denote
it simply\/ $\beta(f)\in P_n$.
\end{theorem}

\begin{proof}
The compact case is trivial, so we assume $M$ non-compact.  We shall
think of $S_n(\bH)_{\Isom^+({\H^n})}$ as being generated by isometry
classes of ideal $n$-simplices.  Triangulate $Z$ and lift the
triangulation to a triangulation of $X$. Note that $Z$ is homeomorphic
to the result of adding a cone on each boundary component of $M_0$.
We may therefore assume that $Z$ is triangulated by first
triangulating $M_0$ and then coning the triangulation at each boundary
component.

As in the proof of
\ref{fund-class}, the map $H_n(\Gamma,\calC)\to \P_n$ can be
identified with the map $H_n(Z)=H_n(C\subdot(X)_\Gamma) \to
H_n(S\subdot(\bH)_{\Isom^+({\H^n})})$ induced by extending the map
$h\colon\calC\to \bH$ to an equivariant map $h'$ of all the vertices
of the triangulation of $X$ to $\bH$.  Given an $n$-simplex of $Z$,
we can lift it to $X$ and then $h'$ maps it to an ideal
$(n+1)$-simplex $\langle z_0,\dots,z_n\rangle\in S_n(\bH)$ which is
well-defined up to the action of $\Gamma$.  In particular it is
well-defined up to isometry. $\mu(\beta_h(f))$ is the sum of these
ideal simplices corresponding to simplices of $Z$.   

Since there are only finitely many $\Gamma$-orbits in $\calC$, it
suffices to see that $\mu(\beta_h(f))$ is unchanged if we change $h$
on just one $\Gamma$-orbit of $\calC$.  Let $v$ be the corresponding
vertex of $Z$ and let $[v,v_1,\dots,v_n]$ be a simplex of $Z$ that
involves $v$.  Let $\langle z,z_1,\dots,z_n\rangle$ be the
corresponding ideal simplex before its change and $\langle
z',z_1,\dots,z_n\rangle$ the ideal simplex after the change. In
$P_n=H_n(S\subdot(\bH)_{\Isom^+(\H^n)})$ the difference of these two
equals
$$\sum_{i=1}^n(-1)^i
\langle{z',z,z_1,\dots,\widetilde{z_i},\dots,z_n}\rangle.$$
But the link of the vertex $v$ is a manifold (it is the corresponding
boundary component of $M_0$) and is triangulated by simplices like
$[v_1,\dots,v_n]$, so the simplex $[v_2,\dots,v_n]$ appears in exactly
two simplices of this link with opposite orientations. Thus the
summand $-\langle z',z,z_2,\dots,z_n\rangle$ of the above sum is
canceled by a corresponding summand from a neighboring simplex. This
is true for all the summands, so the result follows.
\end{proof}

In particular, since volume is well defined on $\P_n$, theorem
\ref{beta of homom} gives a way of defining the ``volume'' of a
homomorphism $f$ as above.  The existence of such a volume in the
3-dimensional case was mentioned in \cite{thurston3}.  In the
three-dimensional case our proof easily gives a little more
information.

\begin{theorem}
Assume $M$ is a hyperbolic $3$-manifold. If $f$ is the homomorphism
corresponding to some Dehn filling $M'$ of $M$ then\/
$\beta(f)=\beta(M')$. If each cusp subgroup of\/ $\Gamma$ has
non-trivial elements $\gamma$ with $f(\gamma)$ parabolic (or trivial)
then\/ $\beta(f)\in\B(\C)$.
\end{theorem}

\begin{proof}
If we are given an ideal triangulation of $M$ then we can compute
$\beta(f)$ as the sum of ideal simplices obtained by taking the
vertices to fixed points of the corresponding cusp subgroups for the
$\Gamma$-action given by $f$. If $f$ corresponds to a Dehn filling
then this gives a degree one ideal triangulation of $M'$ (see e.g.,
\cite{neumann-zagier}) so the sum of the resulting ideal simplices
represents $\beta(M')$. The final sentence follows by the same
argument as in section \ref{proof sect}, except that the isotropy
groups of $\Gamma$ on $\calC$ may map to cyclic rather than  unipotent
subgroups of $\PGL(2,\C)$ and we get a zero map in homology in the
cyclic case because $H_2(\Z)=0$.
\end{proof}

We end this section by discussing to what degree Condition
\ref{condition} is restrictive, and to what extent it can be relaxed.
If we consider $\beta_h(f)\in H_n(\Isom^+(\H^n),\bH)\otimes\Q$ and
$\beta(f)\in\P_n\otimes\Q$ then we need only require that Condition
\ref{condition} is satisfied for a subgroup of finite index in
$\Gamma$, since we can then compute the invariant for this subgroup 
and divide by
the degree of the covering.  In particular, this relaxed condition 
holds for any $f$ if
$n$ is odd and $M$ is a hyperbolic manifold (it suffices that the cusp
groups $P_c$ are virtually polycyclic).  However even without going to
a subgroup of finite index, the condition is not very restrictive.
For instance, if $n=3$ and $M$ is a hyperbolic $3$-manifold then the
cusp groups are isomorphic to $\Z^2$ and the only way the image of
such a group can fail to have a fixed point in $\bH[3]$ is if its 
image is a Klein four-group fixing a point of $\H^3$.

\section{Examples}\label{examples}

In this section, we will look at a few interesting examples that
illustrate our results and conjectures of this paper and
\cite{neumann-yang}. We are grateful to Alan Reid, Frank Calegari, and 
Craig Hodgson for noticing some of these examples.
Many of the calculations in this section are done
with the aid of the software packages Snappea, \cite{weeks}, Snap,
\cite{goodman} and Pari-GP \cite{b-b-c-o}. The number theory
justification of these calculations can be found in the book by Henri
Cohen \cite{cohen}.  

The result that underlies many of our calculations is a theorem of
Borel (reinterpreted in the light of work of Bloch and Suslin --- for
more details, see e.g. \cite{neumann-yang}).  Let $F$ be a number
field and let
$\sigma_1,\overline\sigma_1,\ldots,\sigma_{r_2},\overline\sigma_{r_2}\colon
F\to \C$ be a list of all complex embeddings of $F$.  The ``Borel
regulator'' is the map
$$c_2\colon  \B(F) \longrightarrow  \R^{r_2}$$
defined on generators by
$$c_2([z])=(D_2(\sigma_{1}(z)), \ldots, D_2(\sigma_{r_2}(z))),$$
where $D_2$ is the Bloch-Wigner dilogarithm defined in section
\ref{chern simons}.

\begin{theorem}\label{Borel}
The Borel regulator $c_2$ has kernel the torsion of $\B(F)$ and has
image a maximal sublattice of $\R^{r_2}$.  \qed
\end{theorem}

Thus we can verify computationally whether elements $\alpha$ and
$\beta$ in $\B(F)$ are equal mod torsion by computing whether
$D_2(\sigma(\alpha)) = D_2(\sigma(\beta))$ for all possible embeddings
$\sigma\colon F\subset \C$.  Computing this numerically to sufficient
precision gives absolute proof of the equality of $\alpha$ and $\beta$
modulo torsion if the size of a smallest element of the lattice
$c_2(B(F))\subset \R^{r_2}$ is known. The size of a smallest element can
be bounded in terms of the covolume of the lattice once a rational
basis has been found, and conjectures exist for this covolume which
suggest that one or two digits of precision in calculations would
normally be ample.  These conjectures are open for all but very few
cases, so our examples below cannot be considered to be proved.  But
since we compute to over $50$ digits precision they can probably be
considered to be correct beyond reasonable doubt.

Before we discuss the examples we recall the relationship of the
Chern-Simons invariant with the $\eta$-invariant of
\cite{atiyah-patodi-singer1}.  The $\eta$-invariant is a real valued
invariant defined for compact riemmanian manifolds of dimension
congruent to $3$ modulo $4$. Atiyah, Patodi, and Singer prove in 
\cite{atiyah-patodi-singer2}
\begin{theorem}\label{aps}
$(1/2\pi^2)CS(M)=(3/2)\eta(M)$ $($mod $1/2)$ for any compact riemannian
$3$-manifold $M$.\qed
\end{theorem}
The formula of Theorem \ref{cs-formula} has been improved to a formula
for the eta-invariant $\eta(M)$ in \cite{neumann-meyerhoff} and
\cite{ouyang}, and we use this to give $(3/2)\eta(M)$ rather than
$(1/2\pi^2)CS(M)$ in cases where we have done the computation.  The
normalization $(1/2\pi^2)CS$, which is a well defined invariant modulo
$1/2$, is a fairly standard normalization for $CS$ because for {\em
compact} 3-manifolds $(1/2\pi^2)CS$ is actually well defined modulo
$1$. Nevertheless, our computations of $(1/2\pi^2)CS$ are only valid
modulo $1/2$, even in the compact case.  (The value of
$(1/2\pi^2)CS(M)$ modulo 1 for a compact $3$-manifold is 
determined by $\eta(M)$ and the homology of $M$; see
\cite{atiyah-patodi-singer2}.)

\begin{example}
In Weeks' and Hodgson's census of closed manifolds, included in
Snappea, the one with the smallest volume is $m003(-3,1)$, commonly
known as the Weeks manifold $W$. It is conjectured to be the hyperbolic 
3-manifold of smallest volume. It is obtained via three $(-3,2)$
Dehn surgeries on the alternating 6-crossing link which is a circular
chain of three circles. It can also be obtained as $((5,1),(5,2))$
Dehn surgery on the Whitehead link, which is the description that
Weeks originally found.

The invariant trace field of the Weeks manifold is known to be the
cubic field $k$ of discriminant $-23$ generated by the complex root
with positive imaginary part of the polynomial $x^3 -x + 1$. We denote
this root by $\theta$. It satisfies
$$D_2(\theta)=
0.94270736277692772092129960309221164759032710576688316..., $$
which is the volume of $W$. By Borel's theorem $\B(k)$ is of rank
1. Since $\theta\wedge(1-\theta) = \theta\wedge\theta^3 = 0$,
$[\theta]$ generates $\B(k)\otimes\Q$. The computation of
$D_2(\theta)$ shows that $\beta(W)=[\theta]$ in $\B(k)\otimes\Q$
beyond reasonable doubt.

Since the invariant trace field of $W$ is of odd degree over $\Q$, the
conjecture in the introduction of \cite{neumann-yang} says that
$(1/2\pi^2)CS(W)$ should be irrational.  Computation gives
$$
\frac32\eta(M)=
0.060043066678727155012132615144817756316780200913123686\ldots,
$$
so this is $(1/2\pi^2)CS(W)$ modulo $1/2$.  Using continued fractions
it is easily checked that any rational expression for this would have
to have over $30$ digits in numerator and denominator. This provides
numerical evidence for its irrationality.
\medskip

Alan Reid found that the manifold $M$ obtained via  $(0, 1)$
surgery on the $8_9$-knot has zero Chern-Simons invariant (in fact, 
zero $\eta$-invariant since it admits an orientation reversing symmetry) 
and the
volume of $M$ appears numerically to be 6 times that of $W$. The
Ramakrishnan conjecture (cf.\ \cite{neumann-yang}) suggests that the
Bloch invariant given by $[M]$ be $3([\theta] - [\bar{\theta}])$
modulo torsion and that the invariant trace field of $M$ therefore
contain $\Q(\theta,\overline\theta)$, which is the Galois closure of
$k$.

Computations using Snap show that the invariant trace field of $M$ is
in fact exactly the Galois closure of $k$ and that, at least
numerically, $\beta(M)=3([\theta] - [\bar{\theta}])$ modulo torsion.
\end{example}

\begin{example}
For a number field $F$ with just one complex place there exist
arithmetic 3-manifolds with this field as invariant trace field.  Any
such will give a Bloch invariant which generates $\B(F)\otimes\Q$. For
fields with more than one complex place we do not know how much of
$\B(F)\otimes\Q$ can be generated by Bloch invariants of hyperbolic
3-manifolds.  The following example is of interest in this regard.

The polynomial $f(x)=x^4+x^2-x+1$ is irreducible with roots
$\tau_1^{\pm}=0.54742\ldots\pm 0.58565\ldots i$ and
$\tau_2^{\pm}=-0.54742\ldots\pm 1.12087\ldots i$.  The field
$F=\Q(x)/(f(x))$ is thus of degree $4$ over $\Q$ with two complex
embeddings $\sigma_1,\sigma_2$ up to complex conjugation, one with
image $\sigma_1(F)=\Q(\tau_1^{-})$ and one with image
$\sigma_2(F)=\Q(\tau_2^-)$ (the discriminant of $F$ is $257$, which is
prime, so $F$ has no proper subfield other than $\Q$).  The Bloch
group $\B(F)$ is thus of rank 2 modulo torsion.

It turns out that in the cusped and closed census lists of Snappea
there are a total of five manifolds of volume
$$3.163963228883143983991014715973154484812787671518\ldots,$$
two of them noncompact and three of them compact, and they all have
invariant trace field equal to $\sigma_1(F)$ or its complex
conjugate. (Reversing orientation of a manifold replaces invariant
trace field by the complex conjugate field, so after adjusting
orientations they all have invariant trace field $\sigma_1(F)$.)  Moreover,
there are ten manifolds of volume
$$3.821687586179977739110922224290385516821302495504\ldots,$$
three noncompact and seven compact, and they all have invariant trace
field equal to $\sigma_2(F)$ or its complex conjugate.  One of the seven
compact ones double covers a manifold of half this volume called
$M10(-1,3)$. 

The noncompact manifolds of volume $3.1639\ldots$ are $M032$ and
$M033$ in the 5-simplex census and have 
$(1/2\pi^2)CS(M032)=$ 
$$0.155977016743515161236016645699315761220516234595001\ldots$$
and $(1/2\pi^2)CS(M033)=(1/2\pi^2)CS(M032)-1/4$, while the three
compact ones have $(1/2\pi^2)(CS(M)-CS(M032))$ equal to $-1/3$,
$-3/8$, and $11/60$ respectively.  The noncompact manifolds of volume
$3.82168\ldots$ are $M159$, $M160$, and $M161$ and have 
$(1/2\pi^2)CS(M159)+1/4=(1/2\pi^2)CS(-M160)=(1/2\pi^2)CS(M161)=$
$$ 0.191492799941975387695803880629160641602608213619566\ldots.$$
The compact manifold which double covers $M10(-1,3)$ also has this
Chern Simons invariant.  (In fact $(3/2)\eta(M10(-1,3))=
(.1914927999\ldots-1)/2$.)  The other six compact manifolds of this
volume have Chern-Simons invariants less than $.1914927999\ldots$ by
$1/6,1/3,1/3,1/3,1/3,5/12$ respectively.

Numerical computation shows that the manifolds in the first group all
have the same rational bloch invariant $\beta_1\in\B(F)\otimes\Q$
(more precisely, they all have Bloch invariant $\sigma_1(\beta_1)\in
\B(\sigma_1 F)\otimes\Q$). Similarly the second group gives a class
$\beta_2\in\B(F)\otimes\Q$.  In fact, in the non-compact case the
triangulation gives an exact Bloch invariant in $\B(F)$.  For each of
the two non-compact manifolds of the first class this invariant is
$$\beta_1=2[\frac12(1-\tau^2-\tau^3)] + [1-\tau] 
 + [\frac12(1-\tau^2+\tau^3])\in \B(F)$$
Here $\tau$ denotes the class of $x$ in $F=\Q(x)/(x^4+x^2-x+1)$.
(It is interesting to note that despite the exact equality of Bloch
invariants the Chern-Simons invariants differ by $1/4$ modulo $1/2$.)
Similarly the first two noncompact manifolds of the second group give
the class 
$$\beta_2 = 2[2-\tau-\tau^3]+2[\tau+\tau^2+\tau^3] \in \B(F)$$ 
and the third gives
$$\beta'_2= [\frac14(3+\tau^2)]
+2 [\frac12(\tau^2+\tau^3)] + [\frac14(-3-2\tau-1\tau^2+\tau^3)] +
[\frac1{13}(8-5\tau-2\tau^2-4\tau^3)] \in \B(F).
$$ We do not know if the torsion class
 $\beta_2-\beta'_2$ vanishes.

The Borel regulator map gives:
$$\begin{aligned}
c_2(\beta_1)&= (3.1639632288831439839910147159731544848127876715181,\\
&\qquad -1.4151048972655633406895085877105020361346679596016)\\
c_2(\beta_2)&=(-0.69854408278444071973072661203684276397736670535490,\\
&\qquad 3.8216875861799777391109222242903855168213024955043),
\end{aligned}
$$
proving that $\beta_1$ and  $\beta_2$ generate $\B(F)\otimes\Q$.  

In fact, as we now describe, the whole of $\B(\sigma_1(F))$ is
generated by Bloch invariants of 3-manifolds with invariant trace
field $\sigma_1(F)$.

Searching the closed manifold census for manifolds whose volumes are
small linear combinations of $D_2(\sigma_1(\beta_1)),
D_2(\sigma_1(\beta_2)$ results in five candidates, four of them with
volume $4.396672801932495\ldots$ and one with volume
$5.629382374981847\ldots$. Checking with Snap then confirms that they
all have invariant trace field $\sigma_1(F)$ and their Bloch
invariants in $\B(F)\otimes\Q$ are numerically $(3/2)\beta_1 +
(1/2)\beta_2$ for the four of volume $4.396672801932495\ldots$ and
$2\beta_1+\beta_2$ for the one of volume $5.629382374981847\ldots$.

A similar search for compact manifolds with invariant trace field
$\sigma_2(F)$ yielded no new examples.

The Galois closure $\overline F$ of $F$ is degree $24$ over $\Q$. The
element $\beta_1$ has four distinct Galois conjugates in $\B(\overline
F)\otimes\Q$ (which is of rank 12), and hence in $\B(\C)$. As elements
of $\B(\C)$ these are $\sigma_1(\beta_1),\overline\sigma_1(\beta_1),
\sigma_2(\beta_1),\overline\sigma_2(\beta_1)$.  Their sum is zero (it
follows from Theorem \ref{Borel} that the sum of all Galois conjugates
of any element of $\B(\overline\Q)$ is zero) and they generate a rank
3 subgroup of $\B(\C)$.  Similar remarks apply to $\beta_2$, and one
checks that the two rank 3 subgroups of $\B(\C)$ generated by the
Galois conjugates of $\beta_1$ and $\beta_2$ generate a rank 6
subgroup.  The Bloch invariants of all the above manifolds and their
orientation reversals generate a rank 5 subgroup of this; namely the
subgroup generated by the six elements $\sigma_1(\beta_1)$,
$\overline\sigma_1(\beta_1)$, $\sigma_1(\beta_2),$
$\overline\sigma_1(\beta_2)$, $\sigma_2(\beta_2)$,
$\overline\sigma_2(\beta_2)$, the sum of the last four of which is
zero.
\end{example}

\begin{example}
It is of interest to know to what extent the Bloch invariant
determines Chern Simons invariant in $\R/\pi^2\Z$ rather than in
$\R/\pi^2\Q$.  Let
$$z_1=\frac{3+i-\sqrt{4+2i}}{2},\; z_2=2z_1-2z_1^2+z_1^3/2,\;
z_3=\frac{1+i}{2}.$$ 
Then Snappea shows that the manifolds $M6(1,3)$ and $M11(1,3)$ both
have degree one ideal triangulations using the three simplices with the
above parameters, so their Bloch invariants in $\B(\Q(z_1))$ are equal.
However, their Chern Simons invariants are $11/48$ and $7/48$ which
differ by $1/12$ modulo $1/2$. (In fact $(3/2)\eta(M6(1,3))=-61/48$.) 
This example is especially interesting because these two manifolds are
not commensurable, despite having congruent ideal triangulations.  In
fact, they are both arithmetic with invariant trace field $\Q(i)$, so
by Reid (see e.g., \cite{neumann-reid1}) their commensurability classes
are determined by their invariant quaternion algebras.  But
computation using Snap shows that each of $M6(1,3)$ and $M11(1,3)$ has
quaternion algebra ramified at just two primes, one of which is the
prime dividing $2$ but the other of which divides $5$ or $13$
respectively.   

The volume of these manifolds is
$$1.831931188354438030109207029864768221548298748563344268534.$$
Snappea knows five compact orientable manifolds of this volume and
they lead to several examples like the above. In fact, the five-vertex
graph with edges according to whether the corresponding manifolds have
congruent triangulations with parameters in some quadratic extension
of $\Q(i)$ is a connected graph.  In each case we thus see that the
Bloch invariants of the corresponding manifolds are equal in some
quadratic extension of $\Q(i)$.  We do not know if the Bloch
invariants can be defined in $\B(\Q(i))$, and if so, whether they are
equal there (modulo torsion they are just the element
$2[i]\in\B(\Q(i))$).
\end{example}
\begin{example}
The manifold $X=V3066$ in the seven-simplex cusped census has the
surprising property that the Dehn filled manifolds $X(p,q)$ and
$X(-p,q)$ appear to have equal volume for each $(p,q)$ and to have
Chern-Simons invariants which sum to the apparentally irrational
number  
$$ \alpha=0.02172669391945231711932766534448768004430408\ldots.$$
The manifold $X(1,2)$ has volume 
$$5.137941201873417769841348339474845035649675\ldots$$ and Chern-Simons
invariant $\alpha$.  Moreover, its invariant trace field $k$ is
generated by a root of $x^3+2x-1$ and has discriminant $-59$.  The
Ramakrishnan conjecture would imply that the invariant trace field of
$X(p,q)$  and conjugate invariant trace field of $X(-p,q)$ generate a
field that contains the join of above cubic field and its complex
conjugate, i.e., the Galois closure of this cubic field. 
Experiment suggests that this holds, in fact that the invariant trace
field of  $X(p,q)$ always contains the above cubic field.  For
example, the Ramakrishnan conjecture would imply that $X(-1,2)$ (which
has the same volume as $X(1,2)$ but zero Chern-Simons invariant) must
have invariant trace field containing the Galois closure $K$ of the
above cubic field, and in fact its invariant trace field is exactly
this Galois closure. The Bloch invariant of $X(1,2)$ is in fact
$4(2[\theta]+[1+\theta^2]) \in \B(k)\otimes\Q$, where $\theta$ is the
complex root (with positive imaginary part) of $x^3+2x-1$. And
$X(-1,1)$ and $X(1,1)$ have volume exactly half the above volume,
Chern-Simons invariant $\alpha/2+5/24$ and $\alpha/2-5/24$
respectively, and both also have the above cubic field as invariant
trace field.  The manifold $X$ itself has volume 
$6.2328329776455\ldots$, Chern-Simons invariant $\alpha/2-1/4$, and a
degree 6 invariant trace field of discriminant $-2^659^2$. 
\end{example}

\section{Appendix: Scissors Congruence} \label{appendix}

In this appendix, we will prove that the pre-Bloch group $\P(\C)$
defined in definition \ref{def-bloch} has a more geometric (scissors
congruence) description. Throughout the appendix, by a
face-triangulated polyhedron we mean a convex ideal polyhedron in
$\H^3$ with ideally triangulated faces. As we describe below, we will
allow degenerate (flat) polyhedra, though this is not essential.

One should think of a flat face-triangulated polyhedron as having
infinitesimal thickness, so it is an ideal polygon in $\H^3$ with two
ideal triangulations, one on each ``side''. (A more formal definition
might be to associate a triangulation of the polygon to each of its
two normal directions.) Flat polyhedra occur as follows in
triangulation. In \cite{epstein-penner}, it was shown that a cusped
hyperbolic 3-manifold $M$ can be decomposed into convex ideal
polyhedra. In order to get a triangulation of $M$, one needs to
triangulate each resulting polyhedron. After this triangulation, a
common face of two different polyhedra may now have different
triangulations. In order to make this a true triangulation for $M$,
one needs to insert a flat polyhedron for each such face, which one
should consider to have two sides, with triangulations on each side
to match the triangulations coming from the faces of the two
polyhedra.  One can then triangulate these flat polyhedra into flat 
tetrahedra to complete the triangulation of $M$.  (Of course, by changing 
the triangulations of the polyhedra one may be able to avoid the need for
flat tetrahedra --- it is unknown whether this is always possible.)

In particular, a flat ideal tetrahedron, that is, one with a real
cross ratio parameter $r$, is thus an ideal quadrilateral
triangulated by drawing one diagonal on one ``side'' of it and the
other diagonal on the other ``side''.  To understand which side gets
which diagonal, thicken the flat tetrahedron slightly by deforming $r$
to $r+i\epsilon$ with $\epsilon>0$.

Define a group $\calQ(\C)$ generated by face-triangulated polyhedra
subject to the following relations:\begin{itemize} 
\item for each face-triangulated polyhedron $P$ and isometry $g\in
\Isom^+(\H^3)$, we have $[gP]=[P]$;\item if a face-triangulated
polyhedron $P$ is obtained by gluing two face-triangulated polyhedra
$P_1$ and $P_2$ along a face then $[P] = [P_1]+[P_2]$. The face along
which $P_1$ and $P_2$ are glued together not only should have the same
physical shape, but should also have compatible triangulation.
\end{itemize}

\begin{remark*}{}
There are a couple of points worth noting here:

1.~~The importance of requiring triangular faces in our definition
was made clear to us by a remark of David Kazhdan. Take an ideal
pyramid on an ideal quadrilateral base.  The two ways of cutting the
quadrilateral by a diagonal give two decompositions of the pyramid
into two ideal tetrahedra. If we put these equal then we have made the
flat tetrahedron given by the quadrilateral zero. However the cross
ratios of flat tetrahedra gives $\P(\R)$, the set of all the real
elements in $\P(\C)$. $\P(\R)$ is not trivial in $\P(\C)$. In fact,
after passing to $\B(\C)$, $\B(\R)$ maps onto $\B(\C)_+$. In short,
without requiring triangulated faces, we would be looking at
$\P(\C)/\P(\R)$, which is not what we want.

2.~~It is obvious from the 5-term relation in the definition of
$\P(\C)$ that any real cross ratio can be written as the alternating
sum of complex cross ratios. Thus, if we define a group $\P_0(\C)$ by
replacing $\Z(\C-\{0,1\})$ in Definition \ref{def-bloch} by
$\Z(\C-\R)$, then $\P_0(\C)$ surjects to $\P(\C)$.  It is not hard to
verify that this surjection is an isomorphism (the same holds if $\C$
and $\R$ are replaced by any field and proper subfield).  This
suggests that the group $\calQ(\C)$ remains the same if we use only
non-degenerate polyhedra, which is indeed true and can be proved
without much difficulty. We leave the details to the reader.
\end{remark*}

\begin{proposition}\label{P=Q}
The homomorphism
$$\Phi\colon \P(\C) \rightarrow \calQ(\C)$$
induced by sending $[z]$, for $z=x+iy\in \C-\{0,1\}$ with $y\ge0$, to
any ideal tetrahedron with cross ratio $z$ is well defined and is an
isomorphism.
\end{proposition}

\begin{proof}

First we prove well definedness. 

Given an ideal tetrahedron $\Delta$ and an ordering of its vertices,
the orientation of $\Delta$ induced by that ordering may or may not
agree with the orientation induced on $\Delta$ from $\H^3$ (this
induced orientation makes sense even if $\Delta$ is flat, since we are
giving flat simplices infinitesimal thickness).  We will consider the
tetrahedron $\Delta$ plus the ordering of its vertices to represent
the element $\Delta\in \calQ(\C)$ or $-\Delta\in\calQ(\C)$ according
as these two orientations do or do not agree.  With this convention it
is clear that Equation (\ref{invsim}) in Definition \ref{def-bloch} is
preserved by the map $\Phi$.

Geometrically, equation (\ref{5term}) of Definition \ref{def-bloch}
can be interpreted as follows: given a polyhedron $P_0$ with five
vertices, there is an isometry of $\H^3$ which moves the vertices to
positions $\infty, 0, 1, x, y$, so that the polyhedron can either be
decomposed into three ideal tetrahedra with vertices $\langle\infty,
0, 1, x\rangle$, $\langle\infty, 0, x, y\rangle$ and
$\langle0,1,x,y\rangle$, or it can be decomposed into two ideal
tetrahedra $\langle\infty,0,1,y\rangle$ and $\langle\infty, 1,
x,y\rangle$, and these vertex orderings give the correct orientations
of these five simplices. Then the five term relation (\ref{5term}) of
definition \ref{def-bloch} expresses the equality of these two
decompositions and is thus respected by $\Phi$. It is not hard to
check that if one now permutes the five vertices, the sign of a term
in the five term relation is changed only if the orientation induced
by the vertex ordering of the corresponding simplex has changed, so
the five term relation still expresses the same geometric fact as
before. Therefore the map $\Phi$ is well defined.

If we think of the five-term relation as allowing us to move from one
triangulation of $P_0$ to another, it makes sense to call such a move
a ``cycle move''.

We will now show that $\Phi$ is an isomorphism.  To do so we must show
that any face-triangulated polyhedron $P$ has an ideal triangulation,
that is a subdivision into ideal tetrahedra compatible with the face
triangulations, and moreover, that any two ideal triangulations of $P$
are related by a sequence of cycle moves.

By a {\em triangle of $P$} we shall mean a triangle of the
face-triangulation of a face of $P$. Choose one vertex $v$ of $P$ and
then take the set of cones to $v$ of triangles of $P$ which do not
contain $v$.  These cones are clearly 3-simplices which triangulate
$P$ (there will be flat simplices only if there are triangles of $P$
not containing $v$ in faces that do contain $v$).

Given an arbitrary triangulation of $P$, for each 3-simplex $\Delta$
of the triangulation we can use a cycle move to replace it by the sum
(with appropriate signs or orientations) of the cones to $v$ of the
faces of $\Delta$.  This relates this arbitrary triangulation by cycle
moves to the triangulation just constructed and thus shows that any
two triangulations of $P$ are related by cycle moves, completing the
proof.
\end{proof}

\end{document}